\documentclass[journal]{IEEEtran}
\usepackage{url}

\usepackage[utf8]{inputenc} 
\usepackage[T1]{fontenc}    
\usepackage{booktabs}       
\usepackage{amsfonts}       
\usepackage{nicefrac}       
\usepackage{microtype}      

\usepackage[usenames,dvipsnames]{pstricks}
\usepackage{epsfig}
\usepackage{pst-grad} 
\usepackage{pst-plot} 
\usepackage{amsmath}
\usepackage{amssymb}
\usepackage{theorem}
\usepackage{mathrsfs} 
\usepackage{euscript}
\usepackage{graphicx}
\usepackage{float}
\usepackage{textcomp}
\usepackage{listings}
\usepackage{xcolor}
\usepackage{cases}
\usepackage{algorithm}
\usepackage{multirow}
\usepackage[lofdepth,lotdepth]{subfig}

\definecolor{labelkey}{rgb}{0,0.08,0.45}
\definecolor{refkey}{rgb}{0,0.6,0.0}
\definecolor{Brown}{rgb}{0.45,0.0,0.05}
\definecolor{dgreen}{rgb}{0.00,0.49,0.00}
\definecolor{dblue}{rgb}{0,0.08,0.75}
\definecolor{nido}{rgb}{0.6,0.0,0.4}
\renewcommand{\leq}{\ensuremath{\leqslant}}
\renewcommand{\geq}{\ensuremath{\geqslant}}
\newcommand{\Sum}{\ensuremath{\displaystyle\sum}}

\newcommand{\RPP}{\ensuremath{\left]0,+\infty\right[}}
\newcommand{\HH}{\ensuremath{\mathbb{R}^d}}
\newcommand{\RR}{\ensuremath{\mathbb{R}}}
\newcommand{\NN}{\ensuremath{\mathbb{N}}}

\newcommand{\zeroun}{\ensuremath{\left]0,1\right[}}
\newcommand{\zerounr}{\ensuremath{\left]0,1\right]}}

\newcommand{\minimize}[2]{\ensuremath{\underset{\substack{{#1}}}%
{\text{\rm minimize}}\;\;#2 }}
\newcommand{\pinf}{\ensuremath{{+\infty}}}
\newcommand{\scal}[2]{{\left\langle{{#1}\mid{#2}}\right\rangle}}
\newcommand{\menge}[2]{\big\{{#1}~\big |~{#2}\big\}}
\newcommand{\Menge}[2]{\bigg\{{#1}~\Big|~{#2}\bigg\}}

\newtheorem{theorem}{Theorem}
\newtheorem{assumption}{Assumption}

\newtheorem{example}{Example}
\newtheorem{remark}{Remark}
\newtheorem{problem}{Problem}
\newtheorem{lemma}{Lemma}
\newtheorem{proposition}{Proposition}

\newcommand{\PP}{{\mathsf{P}}}

\begin{document}
\title{Classification and regression using an outer approximation
projection-gradient method}

\author{Michel Barlaud,~\IEEEmembership{Fellow,~IEEE,} 
Wafa Belhajali, 
Patrick L. Combettes~\IEEEmembership{Fellow,~IEEE,} 
and Lionel Fillatre%
\thanks{M. Barlaud, W. Belhajali, and L. Fillatre are with 
Universit\'e C\^ote d'Azur, CNRS, 06900 Sophia Antipolis, 
France.\protect\\%
E-mail: barlaud@i3s.unice.fr, fillatre@i3s.unice.fr}%
\thanks{P. L. Combettes is with the Department of Mathematics, 
North Carolina State University, Raleigh, NC 27695-8205, USA.
\protect\\
E-mail: plc@math.ncsu.edu}}%

\markboth{~}%
{Shell \MakeLowercase{\textit{et al.}}: Bare Demo of IEEEtran.cls
for IEEE Journals}

\maketitle

\begin{abstract}
This paper deals with sparse feature selection and grouping for
classification and regression. The classification or
regression problems under consideration consists in minimizing 
a convex empirical risk
function subject to an $\ell^1$ constraint, a pairwise
$\ell^\infty$ constraint, or a pairwise
$\ell^1$ constraint. Existing work, such as the Lasso
formulation, has focused mainly on Lagrangian penalty
approximations, which often require ad hoc or computationally
expensive procedures to determine the penalization parameter. 
We depart from this approach and
address the constrained problem directly via a splitting method.
The structure of the method is that of the classical 
gradient-projection algorithm, which alternates a 
gradient step on the
objective and a projection step onto the lower level set modeling
the constraint. The novelty of our approach is that the 
projection step is implemented via an outer approximation scheme 
in which the constraint set is approximated by
a sequence of simple convex sets consisting of the intersection of
two half-spaces. Convergence of the iterates generated by the
algorithm is established for a general smooth convex minimization
problem with inequality constraints. Experiments on both synthetic and
biological data show that our method outperforms penalty methods. 
\end{abstract}


\IEEEpeerreviewmaketitle

\section{Introduction}
In many classification and regression problems, the objective is
to select a sparse vector of relevant features.
For example in biological applications, DNA microarray
and new RNA-seq devices provide high dimensional gene expression
(typically  20,000 genes). The challenge is to select the smallest
number of genes (the so-called biomarkers) which are necessary to
achieve accurate biological classification and prediction. A popular
approach to recover sparse feature vectors (under a condition of
mutual incoherence) is to solve a convex optimization problem
involving a data fidelity term $\Phi$ and 
the $\ell^1$ norm \cite{candes,doel,dsUP,tRS}. 
Recent Lasso penalty regularization methods take into account
correlated data using either the pairwise $\ell^1$ penalty
\cite{jpvert,llNC,flasso} or the pairwise $\ell^\infty$ penalty
\cite{oscar} (see also \cite{Figu16} for further developments).
The sparsity or grouping constrained classification problem can be
cast as the minimization of a smooth convex loss subject 
to an $\ell^1$ or a pairwise $\ell^\infty$ constraint, say
$\varphi(w)\leq\eta$. Most of the existing work has focused on 
Lagrangian penalty methods, which aim at solving the unconstrained
problem of minimizing $\Phi+\lambda\varphi$. 
Although, under proper qualification
conditions, there is a formal equivalence between constrained and
unconstrained Lagrangian formulations \cite[Chapter~19]{Livre1},
the exact Lagrange multiplier $\lambda$ can seldom be computed
easily, which leaves the
properties of the resulting solutions loosely defined. 
The main contribution of the present paper is to propose an 
efficient splitting algorithm to solve the constrained formulation
\begin{equation}
\label{e:1}
\minimize{\varphi(w)\leq\eta}{\Phi(w)}
\end{equation}
directly. As discussed in \cite{Imag04}, the bound $\eta$ defining
the constraint can often be determined from prior information on
the type of problem at hand. Our splitting algorithm proceeds by
alternating a gradient step on the smooth classification risk
function $\Phi$ and a projection onto the lower level set
$\menge{w\in\RR^d}{\varphi(w)\leq\eta}$.  The main focus is when
$\varphi$ models the $\ell^1$, pairwise
$\ell^1$ constraint, or pairwise $\ell^\infty$ constraint. 
The projection onto the lower level set is implemented via an outer
projection procedure which consists of successive projections onto
the intersection of two simple half-spaces. The remainder of the
paper is organized as follows.  Section~\ref{sec:classif}
introduces the constrained optimization model.
Section~\ref{sec:pga} presents our new
splitting algorithm, which applies to any constrained smooth
convex optimization problem. In particular we also discuss the
application to regression problems. Section~\ref{sec:exp}
presents experiments on both synthetic and real classical
biological and genomics data base.

\section{Classification risk and convex constraints}
\label{sec:classif}

\subsection{Risk minimization}

We assume that $m$ samples
$(x_i)_{1\leq i\leq m}$ in $\HH$ are available. Typically $m<d$,
where $d$ is the dimension of the feature vector.
Each sample $x_i$ is annotated with a
label $\mathsf{y}_{i}$ taking its value in $\{-1,+1\}$. 
The classification risk associated with a
linear classifier parameterized by a vector $w\in\RR^d$ is
given by
\begin{eqnarray}
\label{wfwEkj7yUj01-20x}
\Phi\colon\RR^d\to\RR\colon w\mapsto\dfrac{1}{m}
\Sum_{i=1}^{m}\phi\big(\mathsf{y}_{i}\scal{x_i}{w}\big).
\end{eqnarray}
We restrict our investigation to convex losses $\phi$ which
satisfy the following assumption.

\begin{assumption}
\label{a:1}
Let $f\colon\RR \to [0,1]$ be an increasing Lipschitz-continuous 
function which is antisymmetric with respect to the point 
$(0,f(0))=(0,{1}/{2})$, integrable, and differentiable at $0$ 
with $f'(0)=\max f'$. 
The loss $\phi\colon\RR\to\RR$ is defined by
\begin{equation}
\label{e:23c4}
(\forall t\in\RR)\quad\phi(t)=-t+\int_{-\infty}^{t}f(s)ds.
\end{equation}
\end{assumption}

The main advantage of this class of smooth losses 
is that it allows us to compute the posterior 
estimation \cite{wnb}.
The function $f$ relates a prediction 
$\scal{x_i}{w}$ of a sample $x_i$ to the posteriori probability
for the class $+1$ via
\begin{eqnarray}
\label{proba}
\widehat{\PP}\big[\mathsf{Y}_i= +1 ~|~ x_i\big]=f(\scal{x_i}{w}).
\end{eqnarray}
This property will be used in Section~\ref{sec:exp} to
compute without any approximation the area under 
the ROC curve (AUC).
The loss $\phi$ in Assumption~\ref{a:1} is convex, everywhere
differentiable with a Lipschitz-continuous derivative, and it is 
twice differentiable at $0$ with $\phi''(0)=\max\phi''$. 
In turn, the function $\Phi$ of \eqref{wfwEkj7yUj01-20x} 
is convex and differentiable, and its gradient
\begin{equation}
\label{e:3}
\nabla\Phi\colon w\mapsto\frac{1}{m}\sum_{i=1}^{m}
f\big(\scal{x_i}{w}\big)x_i
\end{equation}
has Lipschitz constant 
\begin{equation}
\label{c248n3x211-22}
\beta=\frac{f'(0)\sum_{i=1}^{m}\|x_i\|^2}{m}=
\frac{\phi''(0)\sum_{i=1}^{m}\|x_i\|^2}{m}.
\end{equation}
Applications to classification often involve normalized 
features. In this case, \eqref{c248n3x211-22} reduces to
$\beta=f'(0)=\phi''(0)$.
Examples of functions which satisfy Assumption~\ref{a:1} include
that induced by the function $f\colon t\mapsto 1/(1+\exp(-t))$, 
which leads to the logistic loss
$\phi\colon t\mapsto\ln(1+\exp(-t))$,
for which $\phi''(0)=1/4$.
Another example is the Matsusita loss \cite{Matsu}
\begin{eqnarray}
\label{eq:matsuloss}
\phi\colon t\mapsto\frac12\Big(-t+\sqrt{1+t^2}\Big),
\end{eqnarray}
which is induced by $f\colon t\mapsto\big(t/\sqrt{1+t^2}+1\big)/2$.

\subsection{Sparsity model}

In many applications, collecting a sufficient amount of
features to perform prediction is a costly process. The challenge
is therefore to select the smallest number of features (genes or
biomarkers) necessary for an efficient classification and
prediction. The problem can be cast as a constrained optimization 
problem, namely,
\begin{equation}
\label{e:0}
\minimize{\substack{w\in\HH\\ \|w\|_0\leq\delta}}{\Phi(w)},
\end{equation}
where $\|w\|_0$ is the number of nonzero entries of $w$.
Since $\|\cdot\|_0$ is not convex, \eqref{e:0} is
usually intractable and an alternative approach is to use the
norm $\|\cdot\|_1$ as a surrogate, which yields the
Lasso formulation \cite{tRS}
\begin{equation}
\label{e:2}
\minimize{\substack{w\in\HH\\ \|w\|_1\leq\eta}}{\Phi(w)}.
\end{equation}
It has been shown in the context of compressed sensing that under
a so-called restricted isometry property, minimizing with the 
$\|\cdot\|_1$ norm is tantamount to minimizing with the 
$\|\cdot\|_0$ penalty in a sense made precise in \cite{candes}. 

\subsection{Grouping model}

Let us consider the graph $\mathsf{S}$ of connected features 
$(i,j)$. The basic idea is to introduce constraints on 
the coefficients for features
$\omega_i$ and $\omega_j$ connected by an edge in the graph. In
this paper we consider two approaches: directed acyclic graph and
undirected graph. Fused Lasso \cite{flasso}
encourages the coefficients $\omega_i$ and
$\omega_j$ of features $i$ and $j$ connected by an edge in 
the graph to be similar. 
We define the problem of minimizing under the
directed acyclic graph constraint as
\begin{equation}
\label{e:22}
\minimize{\substack{w\in\HH\\ \sum_{(i,j)\in\mathsf{S}} 
{|\omega_i-\omega_j|}\leq\eta}}{\Phi(w)},
\end{equation}
for some suitable parameters $\eta\geq 0$.
In the second, undirected graph, approach \cite{oscar} 
one constrains the coefficients of features $\omega_i$ and
$\omega_j$ connected by an edge using a pairwise $\ell^\infty$
constraint. The problem is to
\begin{equation}
\label{e:222}
\minimize{\substack{w\in\HH\\ \sum_{(i,j)\in\mathsf{S}} 
\max({|\omega_i|, |\omega_j|}) \leq\eta}}{\Phi(w)}.
\end{equation}
To approach the constrained problems \eqref{e:2} and \eqref{e:22},
state of the art methods employ a 
penalized variant \cite{dlSR,dsUP,fht,tRS}. 
In these Lagrangian approaches the objective is to minimize
$\Phi+\lambda\varphi$,
where $\lambda>0$ aims at controlling sparsity and
grouping, and where the constraints are defined by
one of the following (see \eqref{e:2}, \eqref{e:22}, and 
\eqref{e:222})
\begin{equation}
\label{e:def-varphi}
\begin{cases}
\varphi_1=\|\cdot\|_1\\
\varphi_2\colon w\mapsto\sum_{(i,j)\in\mathsf{S}} 
\max({|\omega_i|,|\omega_j|})\\
\varphi_3\colon w\mapsto
\sum_{(i,j)\in\mathsf{S}}{|\omega_i-\omega_j|}.
\end{cases}
\end{equation}
The main drawback of current
penalty formulations resides in the cost associated with the 
reliable computation of the Lagrange multiplier
$\lambda$ using homotopy algorithms \cite{dlSR,fht,hrtzER,myCA}. 
The worst complexity case is $O(3^d)$ \cite{myCA}, which is usually
intractable on real data. 
Although experiments using homotopy algorithms
suggest that the actual complexity is $O(d)$ \cite{myCA}, the
underlying path algorithm remains computationally expensive for 
high-dimensional data sets such as the genomic data set.
To circumvent this computational issue, we propose a new general
algorithm to solve either the sparse ~\eqref{e:2} or the grouping
\eqref{e:22} constrained convex optimization problems directly.

\subsection{Optimization model}
Our classification minimization problem is formally cast as follows.

\begin{problem}
\label{prob:1}
Suppose that $\phi$ satisfies Assumption~\ref{a:1} and let 
$\varphi\colon\RR^d\to\RR$ be convex. Set
\begin{eqnarray}
\label{wfwEkj7yUj01-20xxx}
\Phi\colon\RR^d\to\RR\colon w\mapsto\dfrac{1}{m}
\Sum_{i=1}^{m}\phi\big(\mathsf{y}_{i}\scal{x_i}{w}\big)
\end{eqnarray}
and 
\begin{equation}
\label{e:9}
C=\Menge{w\in\HH}{{\varphi(w)}\leq\eta},
\end{equation}
and let $\beta$ be the Lipschitz constant of $\nabla\Phi$, as 
defined in \eqref{c248n3x211-22}. The problem is to
\begin{equation}
\label{c248n3x211-20}
\minimize{w\in C}{\Phi(w)}.
\end{equation}
\end{problem}

In Section~\ref{sec:exp}, we shall focus on the three instances
of the function $\varphi$ defined in \eqref{e:def-varphi}. We
assume throughout that there exists some $\rho\in\RR$ such that 
$\menge{x\in C}{\Phi(x)\leq\rho}$ is nonempty and bounded, which
guarantees that \eqref{wfwEkj7yUj01-20xxx} has at least one solution. 
In particular, this is true if $\Phi$ or $\varphi$ is coercive. 

\section{Splitting algorithm}
\label{sec:pga}

In this section, we propose an algorithm for solving constrained
classification problem \eqref{c248n3x211-20}.
This algorithm fits in the general category of forward-backward
splitting methods, which have been popular since their
introduction in data processing problem in \cite{Smms05}; see also 
\cite{Siop07,Banf11,Silvia,Bach2011,ML2011}. 
These methods offer flexible
implementations with guaranteed convergence of the sequence of
iterates they generate, a key property to ensure the reliability of
our variational classification scheme.

\subsection{General framework}

As noted in Section~\ref{sec:classif}, $\Phi$ is a 
differentiable convex function
and its gradient has Lipschitz constant $\beta$, where $\beta$ is
given by \eqref{c248n3x211-22}. Likewise, since $\varphi$ is convex
and continuous, $C$ is a closed convex set as a lower level set of
$\varphi$. The principle of a splitting method is to use the
constituents of the problems, here $\Phi$ and $C$, separately.
In the problem at hand, it is natural to use the
projection-gradient method to solve \eqref{c248n3x211-20}. This
method, which is an instance of the proximal forward-backward
algorithm \cite{Smms05}, alternates a gradient step on the
objective $\Phi$ and a projection step onto the constraint set
$C$. It is applicable in the following setting, which captures
Problem~\ref{prob:1}.

\begin{problem}
\label{prob:2}
Let $\Phi\colon\HH\to\RR$ be a differentiable convex function 
such that $\nabla\Phi$ is Lipschitz-continuous with constant
$\beta\in\RPP$, let $\varphi\colon\HH\to\RR$ be a convex function, 
let $\eta\in\RR$, and set $C=\menge{w\in\HH}{\varphi(w)\leq\eta}$. 
The problem is to
\begin{equation}
\label{e:prob2}
\minimize{w\in C}{\Phi(w)}.
\end{equation}
\end{problem}

Let $P_C$ denote the projection operator onto the closed
convex set $C$. 
Given $w_0\in\HH$, a sequence $(\gamma_n)_{n\in\NN}$ of
strictly positive parameters, and a sequence $(a_n)_{n\in\NN}$ in
$\HH$ modeling computational errors in the implementation of the
projection operator $P_C$, the projection-gradient algorithm for
solving Problem~\ref{prob:2} assumes the form
\begin{equation}
\label{c248n3x211-23a}
\begin{array}{l}
\text{for}\;n=0,1,\ldots\\
\left\lfloor
\begin{array}{l}
v_n=w_n-\gamma_n\nabla\Phi(w_n)\\
w_{n+1}=P_C(v_n)+a_n.
\end{array}
\right.\\[2mm]
\end{array}
\end{equation}
We derive at once from \cite[Theorem~3.4(i)]{Smms05} the
following convergence result, which guarantees the convergence of
the iterates.

\begin{theorem}
\label{t:1}
Suppose that Problem~\ref{prob:2} has at least one solution, 
let $w_0\in\HH$, let $(\gamma_n)_{n\in\NN}$ be a sequence in
$\left]0,\pinf\right[$, and let $(a_n)_{n\in\NN}$ be a sequence in
$\HH$ such that
\begin{equation}
\label{wfwEkj7yUj01-20}
\sum_{n\in\NN}\|a_n\|<\pinf,\;\;\inf_{n\in\NN}\gamma_n>0,
\quad\text{and}\quad
\sup_{n\in\NN}\gamma_n<\dfrac{2}{\beta}.
\end{equation}
Then the sequence $(w_n)_{n\in\NN}$ generated by
\eqref{c248n3x211-23a} converges to a solution to
Problem~\ref{prob:2}.
\end{theorem}

Theorem~\ref{t:1} states that the whole sequence of iterates
converges to a solution. Using classical results on the asymptotic
behavior of the projection-gradient method \cite{Lev66b}, we can
complement this result with the following upper bound
on the rate of convergence of the objective value.

\begin{proposition}
\label{r:rate}
In Theorem~\ref{t:1} suppose that $(\forall n\in\NN)$ $a_n=0$. 
Then $\Phi(w_{n+1})-\inf\Phi(C)\leq \vartheta/{n+1}$ for 
some $\vartheta>0$.
\end{proposition}

The implementation of \eqref{c248n3x211-23a} is straightforward
except for the computation of $P_C(v_n)$. Indeed, $C$ is defined in
\eqref{e:9} as the lower level set of a convex function,
and no explicit formula exists for computing the projection onto
such a set in general \cite[Section~29.5]{Livre1}. 
Fortunately, Theorem~\ref{t:1} asserts that 
$P_C(v_n)$ need not be computed exactly. 
Next, we provide an efficient algorithm to compute an approximate
projection onto $C$.

\subsection{Projection onto a lower level set}

Let $p_0\in\HH$, let $\varphi\colon\HH\to\RR$ be a 
convex function, and let $\eta\in\RR$ be such that
\begin{equation}
\label{e:D}
C=\menge{p\in\HH}{\varphi(p)\leq\eta}\neq\varnothing.
\end{equation}
The objective is to compute iteratively the projection $P_C(p_0)$
of $p_0$ onto $C$. The principle of the algorithm is to replace 
this (usually intractable) projection 
by a sequence of projections onto simple outer approximations 
to $C$ consisting of the intersection of two affine half-spaces
\cite{Sico00}.

\begin{figure}
\begin{center}
\scalebox{0.35} 
{
\begin{pspicture}(-14,-7.0)(10.0,10.0)
\definecolor{color2b}{rgb}{0.85,0.85,0.85}
\psellipse[linewidth=0.09,dimen=outer,fillstyle=solid,%
fillcolor=color2b](-4.8,-2.0)(4.8,2.0)
\psellipse[linewidth=0.09,dimen=outer,linestyle=solid,%
linecolor=blue,dash=0.23cm 0.23cm,fillcolor=color2b](-4,-2.0)(8,5.0)
\psline[linewidth=0.07cm,arrowsize=0.30cm 4.01,%
arrowlength=1.4,arrowinset=0.4]{->}(0.72,5.5)(-1.74,1.14)
\psline[linewidth=0.07cm,arrowsize=0.30cm 4.01,linecolor=red,%
arrowlength=1.4,arrowinset=0.4]{->}(0.8,5.5)(1.46,1.76)
\psline[linewidth=0.07cm,arrowsize=0.30cm 4.01,linecolor=blue,%
arrowlength=1.4,arrowinset=0.4]{->}(1.5,1.6)(3.1,4.27)
\psline[linewidth=0.07cm,arrowsize=0.30cm 4.01,linecolor=blue,%
arrowlength=1.4,arrowinset=0.4]{->}(1.5,1.6)(0.44,-0.08)
\psline[linewidth=0.09cm,linestyle=solid,linecolor=red]%
(-13.2,-0.95)(6.5,2.50)
\psline[linewidth=0.09cm,linestyle=dashed,linecolor=blue]%
(-9.2,7.95)(6.5,-1.35)
\psline[linewidth=0.09cm,linecolor=blue](-10.2,6.05)(6.5,-3.85)
\rput(-4.7,-2.0){\Huge 
$C\!=\!\menge{p\in\HH\!}{\!\varphi(p)\leq\eta}$}
\rput(-3.8,-5.2){\Huge\color{blue} $\menge{p\in\HH}{\varphi(p)\leq
\varphi(p_k)}$}
\rput(2.5,2.25){\Huge $p_k$}
\rput(0.6,6.25){\Huge $p_0$}
\rput(0.8,5.6){\Huge $\bullet$}
\rput(1.5,1.6){\Huge $\bullet$}
\rput(-1.77,1.0){\Huge $\bullet$}
\rput(-1.81,0.3){\Huge $p_{k+1}$}
\rput(0.39,-0.22){\Huge $\bullet$}
\rput(1.93,-0.22){\Huge $p_{k+1/2}$}
\rput(3.6,4.6){\Huge\color{blue} $s_k$}
\rput(-10.2,3.8){\Huge\color{blue} $H(p_k,p_{k+1/2})$}
\rput(4.9,1.35){\Huge\color{red} $H(p_0,p_k)$}
\end{pspicture}
}
\caption{A generic iteration for
computing the projection of $p_0$ onto $C$. At iteration $k$, the
current iterate is $p_k$ and $C$ is contained in the half-space
$H(p_0,p_{k})$, onto which $p_k$ is the projection of $p_0$ (see
\eqref{eq:halfspaceH}). If $\varphi(p_k)>\eta$, a subgradient
vector $s_k\in\partial\varphi(p_k)$ is in the normal cone 
to the lower level set 
$\menge{p\in\HH}{\varphi(p)\leq\varphi(p_k)}$ at $p_k$, 
and the subgradient projection $p_{k+1/2}$ of $p_k$ is defined by
\eqref{Vge84Hy03-04a}; it is the projection of $p_k$
onto the half-space $H(p_k,p_{k+1/2})$ which contains $C$.
The update $p_{k+1}$ is the projection of $p_0$ onto
$H(p_0,p_{k})\cap H(p_k,p_{k+1/2})$.}
\label{fig:haugazeau}
\end{center}
\end{figure}
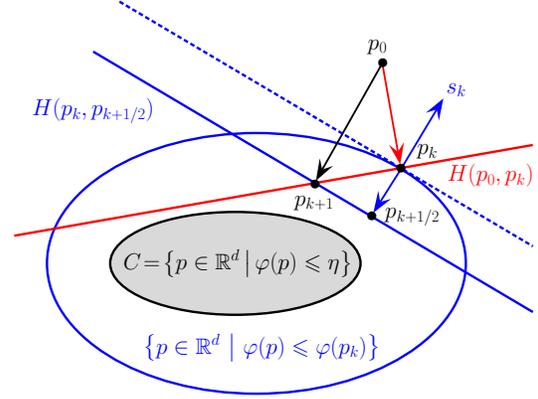

We first recall that $s\in\HH$ is called a subgradient of
$\varphi$ at $p\in\HH$ if \cite[Chapter~16]{Livre1}
\begin{equation}
(\forall y\in\HH)\quad\scal{y-p}{s}+\varphi(p)\leq\varphi(y).
\end{equation}
The set of all subgradients of $\varphi$ at $p$ is denoted by
$\partial\varphi(p)$. If $\varphi$ is differentiable at $p$, this
set reduces to a single vector, namely the gradient 
$\nabla\varphi(p)$.
The projection $P_C(p_0)$ of 
$p_0$ onto $C$ is characterized by 
\begin{equation}
\label{eq:condition-projection}
\begin{cases}
P_C(p_0)\in C\\
(\forall p\in C)\;\scal{p-P_C(p_0)}{p_0-P_C(p_0)}\leq 0.
\end{cases}
\end{equation}
Given $x$ and $y$ in $\HH$, define a closed affine half-space 
$H(x,y)$ by
\begin{equation}
\label{eq:halfspaceH}
H(x,y)=\menge{p\in\HH}{\scal{p-y}{x-y}\leq 0}.
\end{equation}
Note that $H(x,x)=\RR^d$ and, if $x\neq y$, $H(x,y)$ is the 
closed affine half-space onto which the projection of $x$ is $y$.
According to \eqref{eq:condition-projection}, 
$C\subset H(p_0,P_C(p_0))$.
The principle of the algorithm is as follows
(see Fig.~\ref{fig:haugazeau}). At iteration $k$, if
$\varphi(p_k)\leq\eta$, then $p_k\in C$ and the algorithm 
terminates with $p_k=P_C(p_0)$. 
Indeed, since $C\subset H(p_0,p_k)$ \cite[Section~5.2]{Sico00}
and $p_k$ is the projection
of $p_0$ onto $H(p_0,p_k)$, we have
$\|p_0-p_k\|\leq\|p_0-P_C(p_0)\|$. 
Hence $p_k\in C$
$\Leftrightarrow$ $p_k=P_C(p_0)$, i.e., $\varphi(p_k)\leq\eta$
$\Leftrightarrow$ $p_k=P_C(p_0)$. 
Otherwise, one first computes the 
so-called subgradient projection of $p_k$ onto $C$. 
Recall that, given $s_k\in\partial\varphi(p_k)$, the subgradient
projection of $p_k$ onto $C$ is \cite{Livre1,Cens82,Imag97}
\begin{equation}
\label{Vge84Hy03-04a}
p_{k+1/2}=
\begin{cases}
p_k+\displaystyle{\frac{\eta-\varphi(p_k)}{\|s_k\|^2}}s_k
&\text{if}\;\;\varphi(p_k)>\eta\\
p_k&\text{if}\;\;\varphi(p_k)\leq\eta.
\end{cases}
\end{equation}
As noted in \cite{Imag97}, the closed half-space 
$H(p_k,p_{k+1/2})$ serves as an outer approximation to $C$ at 
iteration $k$, i.e., $C\subset H(p_k,p_{k+1/2})$; moreover
$p_k\notin C$ $\Rightarrow$ $p_k\notin H(p_k,p_{k+1/2})$.
Thus, since we have also seen that $C\subset H(p_0,p_k)$, we have
\begin{equation}
C\subset C_k,\quad\text{where}\quad
C_k=H(p_0,p_k)\cap H(p_k,p_{k+1/2}). 
\end{equation}
The update $p_{k+1}$ is computed as the projection 
of $p_0$ onto the outer approximation $C_k$. 
As the following lemma from~\cite{Haug68}
(see also \cite[Corollary~29.25]{Livre1}) shows, this
computation is straightforward.

\begin{lemma}
\label{l:haugazeauy}
Let $x$, $y$, and $z$ be points in $\HH$ such that
\begin{equation}
H(x,y)\cap H(y,z)\neq\varnothing.
\end{equation}
Moreover, set $a=x-y$, $b=y-z$, $\chi=\scal{a}{b}$,
$\mu=\|a\|^2$, $\nu=\|b\|^2$, and
$\rho=\mu\nu-\chi^2$.
Then the projection of $x$ onto $H(x,y)\cap H(y,z)$ is
\begin{equation}
\label{e:2010-09-07v}
Q(x,y,z)=
\begin{cases}
z\!\!&\:\text{if}\;\rho=0\;\text{and}\;\chi\!\geq\! 0\\[+3mm]
\displaystyle x-\bigg(1+\frac{\chi}{\nu}\bigg)b
\!\!&\:\text{if}\;\rho\!>\!0\;\text{and}\;
\chi\nu\geq\rho\\[+4mm]
\displaystyle y+\frac{\nu}{\rho}
\big(\chi a-\mu b\big) 
\!\!&\:\text{if}\;\rho\!>\!0\;\text{and}\;\chi\nu\!<\!\rho.
\end{cases}
\end{equation}
\end{lemma}

To sum up, the projection of $p_0$ onto the set $C$ of \eqref{e:D} 
will be performed by executing the following routine.
\begin{equation}
\label{e:24}
\begin{array}{l}
\text{for}\;k=0,1,\ldots\\
\left\lfloor
\begin{array}{l}
\text{if}\;\varphi(p_k)\leq\eta\\
\!\!\begin{array}{l}
\left\lfloor
\text{terminate.}
\right.\\
\end{array}\\
\zeta_k=\eta-\varphi(p_k)\\
s_k\in\partial\varphi(p_k)\\
p_{k+1/2}=p_k+\zeta_ks_k/\|s_k\|^2\\
p_{k+1}=Q(p_0,p_k,p_{k+1/2}).
\end{array}
\right.\\[2mm]
\end{array}
\end{equation}
The next result from 
\cite[Section~6.5]{Sico00} guarantees the convergence of the 
sequence $(p_k)_{k\in\NN}$ generated by \eqref{e:24} to the 
desired point.

\begin{proposition}
\label{p:22}
Let $p_0\in\HH$, let $\varphi\colon\HH\to\RR$ be a convex
function, and let $\eta\in\RR$ be such that
$C=\menge{p\in\HH}{\varphi(p)\leq\eta}\neq\varnothing$.
Then either \eqref{e:24} terminates in a finite
number of iterations at $P_C(p_0)$ or it generates an 
infinite sequence $(p_k)_{k\in\NN}$ such that $p_k\to P_C(p_0)$.
\end{proposition}

\begin{figure}
\begin{center}
\includegraphics[width=\linewidth]{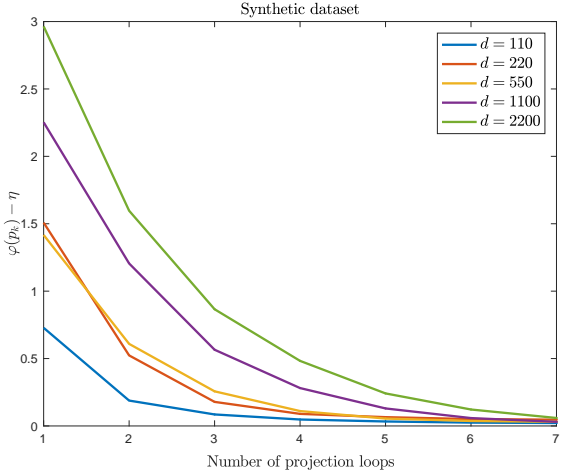}
\end{center}
\caption{Convergence of the projection loop (7 iterations).}
\label{fig:inner}
\end{figure}

To obtain an implementable version of the conceptual algorithm
\eqref{c248n3x211-23a}, consider its $n$th iteration and the
computation of the approximate projection $w_{n+1}$ of $v_n$ onto
$C$ using \eqref{e:24}. We first initialize \eqref{e:24} with 
$p_0=v_n$, and then execute only $K_n$ iterations of it. 
In doing so, we approximate the exact projection onto $C$ by 
the projection $p_{K_n}$
onto $C_{K_n-1}$. The resulting error
is $a_n=P_C(p_0)-p_{K_n}$. According to Theorem~\ref{t:1}, 
this error must be controlled so as to yield overall a summable 
process. First, since $P_C$ is nonexpansive 
\cite[Proposition~4.16]{Livre1}, we have 
\begin{equation}
\|P_C(p_0)-P_C(p_{K_n})\|\leq\|p_0-p_{K_n}\|\to 0.
\end{equation}
Now suppose that $\varphi(p_{K_n})>\eta$ (otherwise we are done).
By convexity, $\varphi$ is Lipschitz-continuous relative to
compact sets \cite[Corollary~8.41]{Livre1}.
Therefore there exists $\zeta>0$ such that
$0<\varphi(p_{K_n})-\eta=
\varphi(p_{K_n})-\varphi(P_C(p_0))\leq\zeta\|p_{K_n}-P_C(p_0)\|
\to 0$.
In addition, assuming that $\text{int}(C)\neq\varnothing$,
using standard error bounds on convex inequalities
\cite{Lewi98}, there exists a constant $\xi>0$ such that
\begin{equation}
\|p_{K_n}-P_C(p_{K_n})\|\leq\xi\big(\varphi(p_{K_n})-\eta\big)
\to 0.
\end{equation}
Thus, 
\begin{align}
\|a_n\|
&=\|P_C(p_0)-p_{K_n}\|\nonumber\\
&\leq\|P_C(p_0)-P_C(p_{K_n})\|+\|P_C(p_{K_n})-p_{K_n}\|\nonumber\\
&\leq\|p_0-p_{K_n}\|+\xi\big(\varphi(p_{K_n})-\eta\big).
\end{align}
Thus, is suffices to take $K_n$ large enough so that, for
instance, we have
$\|p_0-p_{K_n}\|\leq\xi_1/n^{1+\epsilon}$
and $\varphi(p_{K_n})-\eta\leq\xi_2/n^{1+\epsilon}$ for some
$\xi_1>0$, $\xi_2>0$, and $\epsilon>0$. 
This will guarantee that $\sum_{n\in\NN}\|a_n\|<\pinf$
and therefore, by Theorem~\ref{t:1}, the convergence of the 
sequence $(w_n)_{n\in\NN}$ generated by the following algorithm to 
a solution to Problem~\ref{prob:2}.
\begin{equation}
\label{wfwEkj7yUj01-23x}
\begin{array}{l}
\text{for}\;n=0,1,\ldots\\
\left\lfloor
\begin{array}{l}
v_n=w_n-\gamma_n\nabla\Phi(w_n)\\
p_0=v_n\\
\begin{array}{l}
\hskip -2mm\text{for}\;k=0,1,\ldots,K_n-1\\
\left\lfloor
\begin{array}{l}
\zeta_k=\eta-\varphi(p_k)\\
\text{if}\;\zeta_k\geq 0\\
\begin{array}{l}
\left\lfloor
\text{terminate.}
\right.\\
\end{array}\\
s_k\in\partial\varphi(p_k)\\
p_{k+1/2}=p_k+\zeta_ks_k/\|s_k\|^2\\
p_{k+1}=Q(p_0,p_k,p_{k+1/2})
\end{array}
\right.\\[2mm]
\end{array}\\
w_{n+1}=p_{K_n}.
\end{array}
\right.\\[2mm]
\end{array}
\end{equation}
Let us observe that, from a practical standpoint, we have found 
the above error analysis not to be required in our experiments
since an almost exact projection is actually obtainable with a few
iterations of \eqref{e:24}. 
For instance, numerical simulations 
(see Fig.~\ref{fig:inner}) on the synthetic
data set described in Section~\ref{subsec:synthetic} show 
that \eqref{e:24} yields in about 
$K_n\approx 7$ iterations a point very close to the exact 
projection of $p_0$ onto $C$. 
Note that the number of iterations of \eqref{e:24} does 
not depend on the dimension $d$. 

\begin{remark}[multiple constraints]
{\rm We have presented above the case of a single constraint,
since it is the setting employed in subsequent sections. 
However, the results of \cite[Section~6.5]{Sico00} enable 
us to extend this approach to problems with $p$ constraints, see
Appendix~A.}
\end{remark}

\subsection{Application to Problem~\ref{prob:1}}
It follows from \eqref{e:3}, \eqref{e:2010-09-07v}, 
and \eqref{e:24}, that \eqref{wfwEkj7yUj01-23x} for the 
classification problem can be written explicitly as follows,
where $\varepsilon$ is an arbitrarily small number in $\zeroun$
and where $\beta$ is given by \eqref{c248n3x211-22}.
\begin{equation}
\label{wfwEkj7yUj01-23}
\begin{array}{l}
\text{for}\;n=0,1,\ldots\\
\left\lfloor
\hskip -1mm\begin{array}{l}
\gamma_n\in[\varepsilon,(2-\varepsilon)/\beta]\\
v_n=w_n-\dfrac{\gamma_n}{m}\Sum_{i=1}^{m}\mathsf{y}_{i}\phi'
\big(\mathsf{y}_{i}\scal{x_i}{w_n}\big)x_i\\[-1mm]
p_0=v_n\\
\begin{array}{l}
\hskip -2mm\text{for}\;k=0,1,\ldots,K_n-1\\
\hskip -1mm\left\lfloor
\hskip -1mm\begin{array}{l}
\zeta_k=\eta-\varphi(p_k)\\
\text{if}\;\zeta_k\geq 0\\
\!\!\begin{array}{l}
\left\lfloor
\text{terminate.}
\right.\\
\end{array}\\
s_k\in\partial\varphi(p_k)\\
p_{k+1/2}=p_k+\zeta_ks_k/\|s_k\|^2\\
\chi_k=\scal{p_0-p_k}{p_k-p_{k+1/2}}\\
\mu_k=\|p_0-p_k\|^2\\
\nu_k=\|p_k-p_{k+1/2}\|^2\\
\rho_k=\mu_k\nu_k-\chi_k^2\\
\text{if}\;\rho_k=0\;\text{and}\;\chi_k\geq 0\\
\!\!\begin{array}{l}
\left\lfloor
p_{k+1}=p_{k+1/2}
\right.\\[2mm]
\end{array}\\
\text{if}\;\rho_k>0\;\text{and}\;\chi_k\nu_k\geq\rho_k\\
\!\!\begin{array}{l}
\left\lfloor
p_{k+1}=\displaystyle p_0+
\bigg(1+\frac{\chi_k}{\nu_k}\bigg)\big(p_{k+1/2}-p_k\big)
\right.\\[2mm]
\end{array}\\[4mm]
\text{if}\;\rho_k>0\;\text{and}\;\chi_k\nu_k<\rho_k\\
\!\!\begin{array}{l}
\left\lfloor
\!\!\begin{array}{l}
p_{k+1}\!=\!\displaystyle p_k\!+\!\frac{\nu_k}{\rho_k}
\Big(\chi_k\big(p_0\!-\!p_k\big)\!+\!\mu_k\big(p_{k+1/2}\!-\!
p_k\big)\Big)
\end{array}
\right.\\[2mm]
\end{array}
\end{array}
\right.\\[2mm]
\end{array}\\
w_{n+1}=p_{K_n}.
\end{array}
\right.\\[2mm]
\end{array}
\end{equation}
A subgradient 
of $\varphi_1$ at $(\xi_i)_{1\leq i\leq d}\in\RR^d$ is
$s=(\operatorname{sign}(\xi_i))_{1\leq i\leq d}$, where 
\begin{equation}
\label{e:s7}
\operatorname{sign}\colon\xi\mapsto
\begin{cases}
1&\text{if}\;\;\xi>0\\
0&\text{if}\;\;\xi=0\\
-1&\text{if}\;\;\xi<0.
\end{cases}
\end{equation}
The $i$th
component of a subgradient of $\varphi_2$ at
$(\xi_i)_{1\leq i\leq d}\in\RR^d$ is given by
\begin{equation}
\label{e:s77}
\Sum_{(i,j)\in\mathsf{S}} 
\begin{cases}
\operatorname{sign}(\xi_i)&\text{if}\;\;|\xi_i|\geq|\xi_j|\\
0&\text{otherwise}.
\end{cases}
\end{equation}
The $i$th
component of a subgradient of $\varphi_3$ at
$(\xi_i)_{1\leq i\leq d}\in\RR^d$ is given by
\begin{equation}
\label{e:s777}
\Sum_{(i,j)\in\mathsf{S}}
\begin{cases}
\operatorname{sign}(\xi_i-\xi_j)
&\text{if}\;\;\xi_i\neq\xi_j\\
0&\text{otherwise}.
\end{cases}
\end{equation}

\subsection{Application to regression}
A common approach in regression is to learn $w\in\HH$ by
employing the quadratic loss 
\begin{equation}
\label{e:regression}
\Psi\colon\RR^d\to\RR\colon w\mapsto\dfrac{1}{2m}
\Sum_{i=1}^{m}\big|\scal{x_i}{w}-\mathsf{y}_{i}\big|^2
\end{equation}
instead of the function $\Phi$ of \eqref{wfwEkj7yUj01-20xxx} in
Problem~\ref{prob:2}. Since $\Psi$ is convex and has a
Lipschitz-continuous gradient with constant $\beta=\sigma_1^2$,  
where $\sigma_1$ is the largest singular value of the matrix 
$[x_1|\cdots|x_m]$, it suffices to change the definition
of $v_n$ in \eqref{wfwEkj7yUj01-23} by
\begin{equation}
\label{wfwEkj7yUj01-24}
v_n=w_n-\dfrac{\gamma_n}{m}
\Sum_{i=1}^{m}(\scal{x_i}{w_n}-\mathsf{y}_{i}\big)x_i.
\end{equation}

\section{Experimental evaluation}
\label{sec:exp}
We illustrate the performance of the proposed constrained 
splitting method on both synthetic and real data sets.

\subsection{Synthetic data set }
\label{subsec:synthetic}

We first simulate a simple regulatory network in genomic
described in \cite{llNC}. A genomic
network is composed of regulators (transcription factors,
cytokines, kinase, growth factors, etc.) and the genes they
regulate. Our notation is as follows:
\begin{itemize}
\item 
$m$: number of samples.
\item 
$N_\text{reg}$: number of regulators.
\item 
$N_\text{g}$: number of genes per regulator.
\item 
$d=N_\text{reg}(N_\text{g}+1)$.
\end{itemize}
The entry $\xi_{i,j}$ of the matrix  $X=[x_1|\cdots|x_m]^\top$, 
composed of $m$
rows and $d$ columns, is as follows.
\begin{enumerate}
\item 
The $r$th regulator of the $i$th sample is
\begin{eqnarray*}
\xi_{i,\text{reg}_r}
=\xi_{i,N_\text{g}(r-1)+r}
=\xi_{i,r(N_\text{g}+1)-N_\text{g}}\sim \mathcal{N}(0,1).
\end{eqnarray*}
This defines $\xi_{i,j}$ for $j$ of the form
$r(N_\text{g}+1)-N_\text{g}$.
\item
The genes associated with  $\xi_{i,\text{reg}_r}$ have a 
joint bivariate normal distribution with a correlation of $\rho=0.7$
\begin{eqnarray*}
\xi_{i,r(N_\text{g}+1)-N_\text{g}+k}\sim
\mathcal{N}\big(\varrho\,\xi_{i,\text{reg}_r},1-\varrho^2\big).
\end{eqnarray*}
This defines $\xi_{i,j}\neq r(N_\text{g}+1)-N_\text{g}$.
\end{enumerate}
The regression response $\mathsf{Y}$ is given by 
$\mathsf{Y}=Xw+\varepsilon$, where $\varepsilon \sim
\mathcal{N}(0,\sigma^2)$ with $\sigma=2$. 

\begin{example}
\label{ex:1} 
\rm
In this example, we consider that $9$ genes regulated by the same 
regulators are activated and $1$ gene is 
inhibited. The true regressor is defined as
\[
\begin{array}{l}
w=\Bigg(5, \underbrace {\frac{5}{\sqrt{10}}, \ldots }_{9},
\underbrace{\frac{-5}{\sqrt{10}}, \ldots }_{1}  ,   -5, 
\underbrace
{\frac{-5}{\sqrt{10}}, \ldots }_{9},\underbrace{\frac{5}{\sqrt{10}}
, \ldots }_{1} , 3,\\
\underbrace {\frac{3}{\sqrt{10}}, \ldots}_{9},
\underbrace{\frac{-3}{\sqrt{10}}  , \ldots }_{1}  ,-3, \underbrace
{\frac{-3}{\sqrt{10}}, \ldots }_{9},
\underbrace{\frac{3}{\sqrt{10}}, \ldots }_{1}, 0, \ldots\!, 0
\Bigg).
\end{array}
\]
\end{example}

\begin{example}
\label{ex:2}  
\rm
We consider that $8$ genes regulated by the
same regulators are activated and $2$ genes are inhibited.
The true regressor is defined as
\[
\begin{array}{l}
w=\Bigg(5, \underbrace {\frac{5}{\sqrt{10}}, \ldots }_{8},
\underbrace{\frac{-5}{\sqrt{10}}, \ldots }_{2}  ,   -5,
\underbrace {\frac{-5}{\sqrt{10}}, \ldots
}_{8},\underbrace{\frac{5}{\sqrt{10}} , \ldots }_{2} , 3,\\
\underbrace {\frac{3}{\sqrt{10}}, \ldots}_{8},
\underbrace{\frac{-3}{\sqrt{10}}  , \ldots }_{2}  ,-3,
\underbrace {\frac{-3}{\sqrt{10}}, \ldots }_{8},
\underbrace{\frac{3}{\sqrt{10}}, \ldots }_{2}, 0, \ldots\!, 0\Bigg).
\end{array}
\]
\end{example}

\begin{example}
\label{ex:3} 
\rm
This example is similar to Example~\ref{ex:1}, but we consider that
$7$ genes regulated by the same regulators are activated and $3$
genes are inhibited. The true regressor is 
\[\begin{array}{l}
w=\Bigg(5, \underbrace {\frac{5}{\sqrt{10}}, \ldots }_{7},
\underbrace{\frac{-5}{\sqrt{10}}, \ldots }_{3}  ,   -5, \underbrace
{\frac{-5}{\sqrt{10}}, \ldots }_{7},\underbrace{\frac{5}{\sqrt{10}}, 
\ldots }_{3} , 3, \\
\underbrace {\frac{3}{\sqrt{10}}, \ldots}_{7},
\underbrace{\frac{-3}{\sqrt{10}}  , \ldots }_{3}  ,-3, \underbrace
{\frac{-3}{\sqrt{10}}, \ldots }_{7},
\underbrace{\frac{3}{\sqrt{10}}, \ldots }_{3}, 0, \ldots\!, 0\Bigg).
\end{array}
\]
\end{example}

\subsection{Breast cancer data set}
\label{subsec:breast}
We use the breast
cancer data set \cite{van}, which consists of gene expression
data for 8,141 genes in 295 breast cancer tumors (78 metastatic
and 217 non-metastatic).  In the time comparison evaluation, we 
select a subset of the 8141 genes (range 3000 to 7000) using a
threshold on the mean of the genes. We use the network provided in
\cite{network} with $p=639$ pathways as graph constraints in our
classifier. In biological applications, pathways are genes
grouped according to their biological functions \cite{network,llNC}. 
Two genes are connected if they belong to the same pathway. 
Let $\mathsf{S_i}$ be the subset of genes that are connected to 
gene $i$. In this case, we have a subset of only 40,000 connected
genes in $\mathsf{S_i}$. 
Note that we compute the subgradient \eqref{e:s77} only on the
subset $\mathsf{S_i}$ of connected genes. 

\subsection{Comparison between penalty method and our
$\ell^1$ constrained method for classification}

\begin{figure}[H]
\begin{center}
\includegraphics[width=\linewidth]{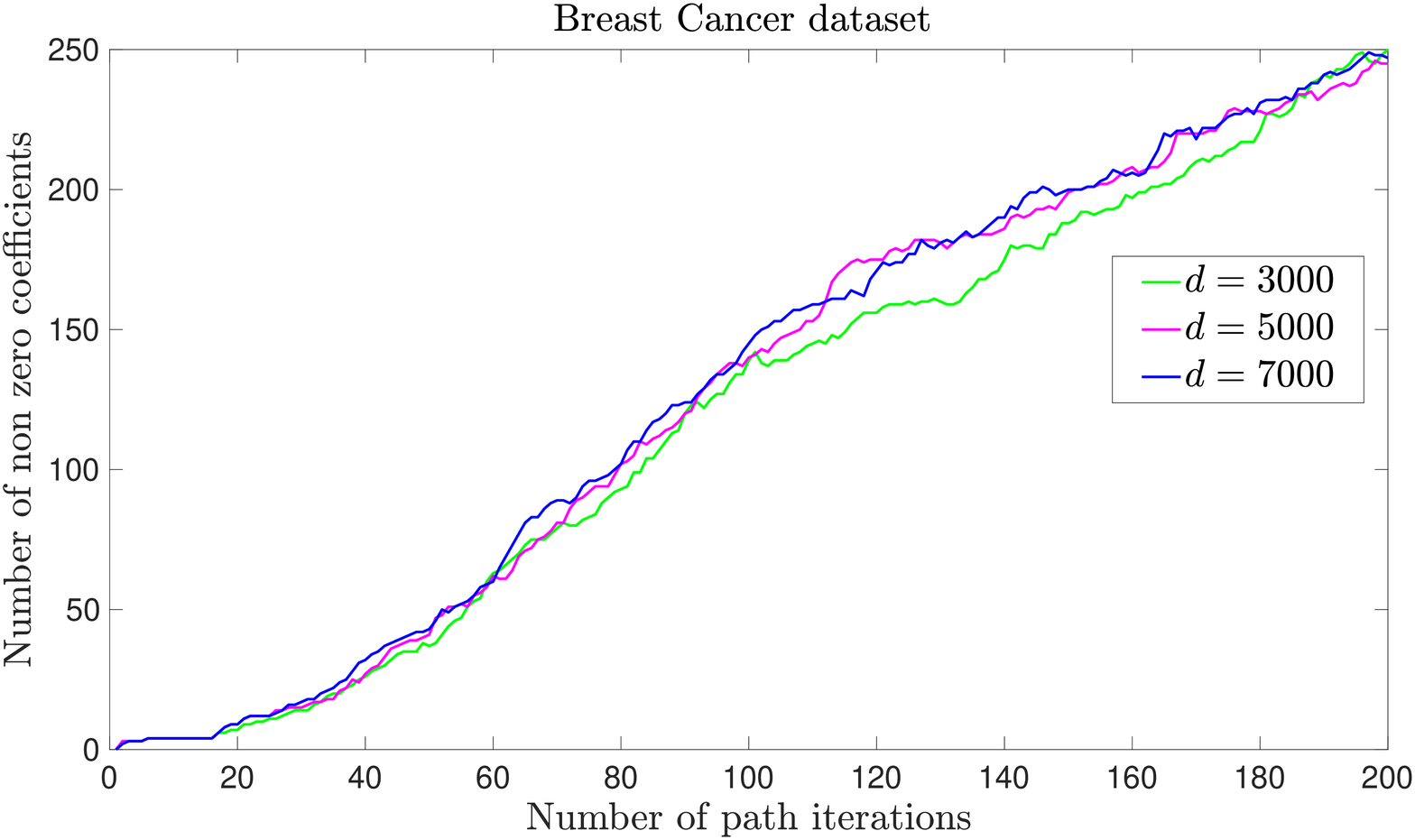}
\end{center}
\caption{Glmnet: Number of nonzero coefficients as a function of 
$n_{\lambda}$.}
\label{fig:lambda}
\end{figure}

First, we compare with the penalty approach using glmnet
MATLAB software \cite{glmnet} on the breast cancer data set described
in Section~\ref{subsec:breast}. We tuned the number of path
iterations $n_{\lambda}$ for glmnet for different values of the
feature dimension. The number of nonzero coefficients $\|w\|_0$
increases with $n_{\lambda}$.
The glmnet method requires typically 200 path iterations
or more (see Fig.~\ref{fig:lambda}). 

\begin{figure}[H]
\begin{center}
\includegraphics[width=\linewidth]{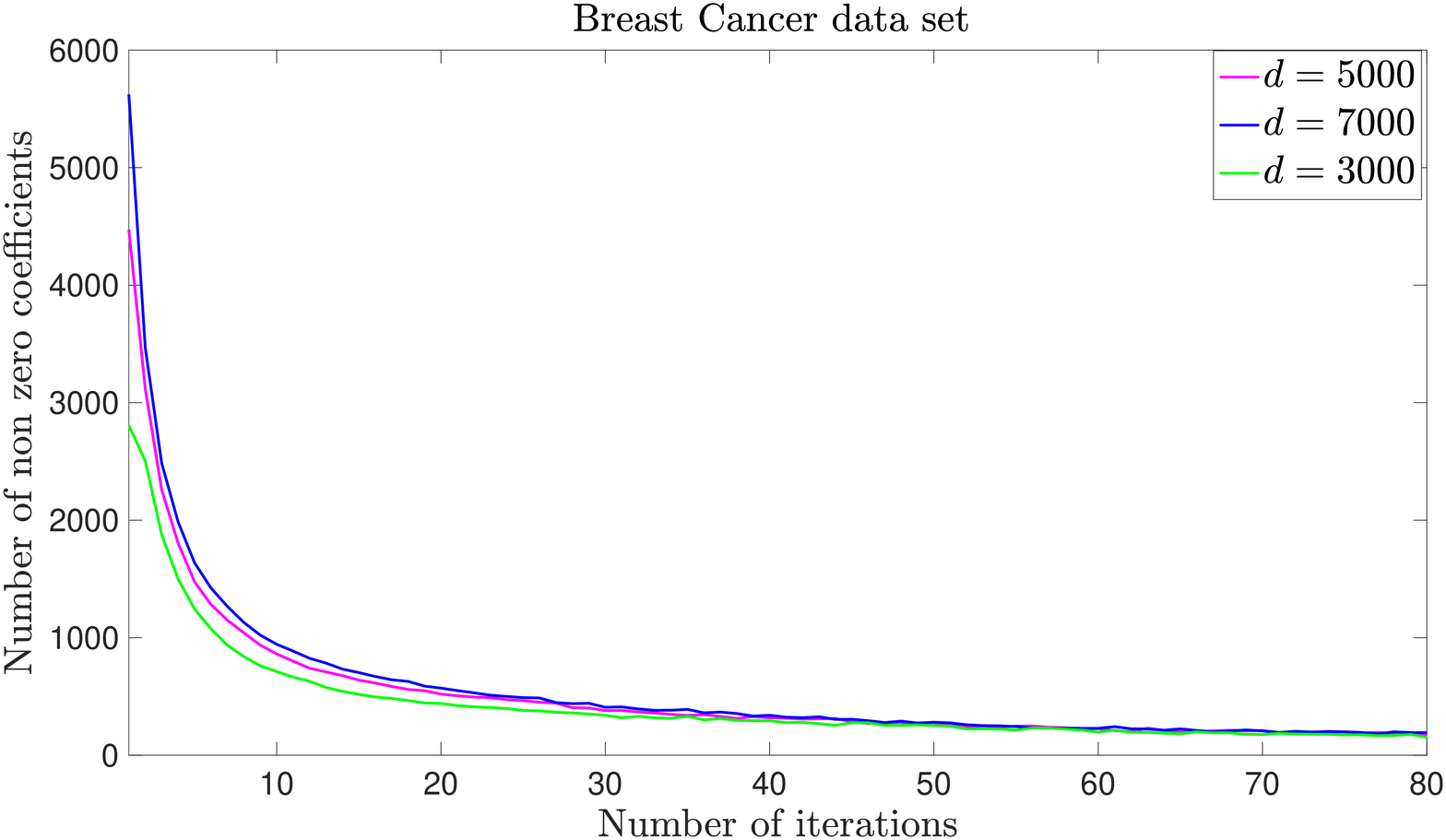}
\end{center}
\caption{Number of nonzero coefficients as a function of the
number of iterations.}
\label{fig:sparse}
\end{figure}

\begin{figure}
\begin{center}
\includegraphics[width=\linewidth]{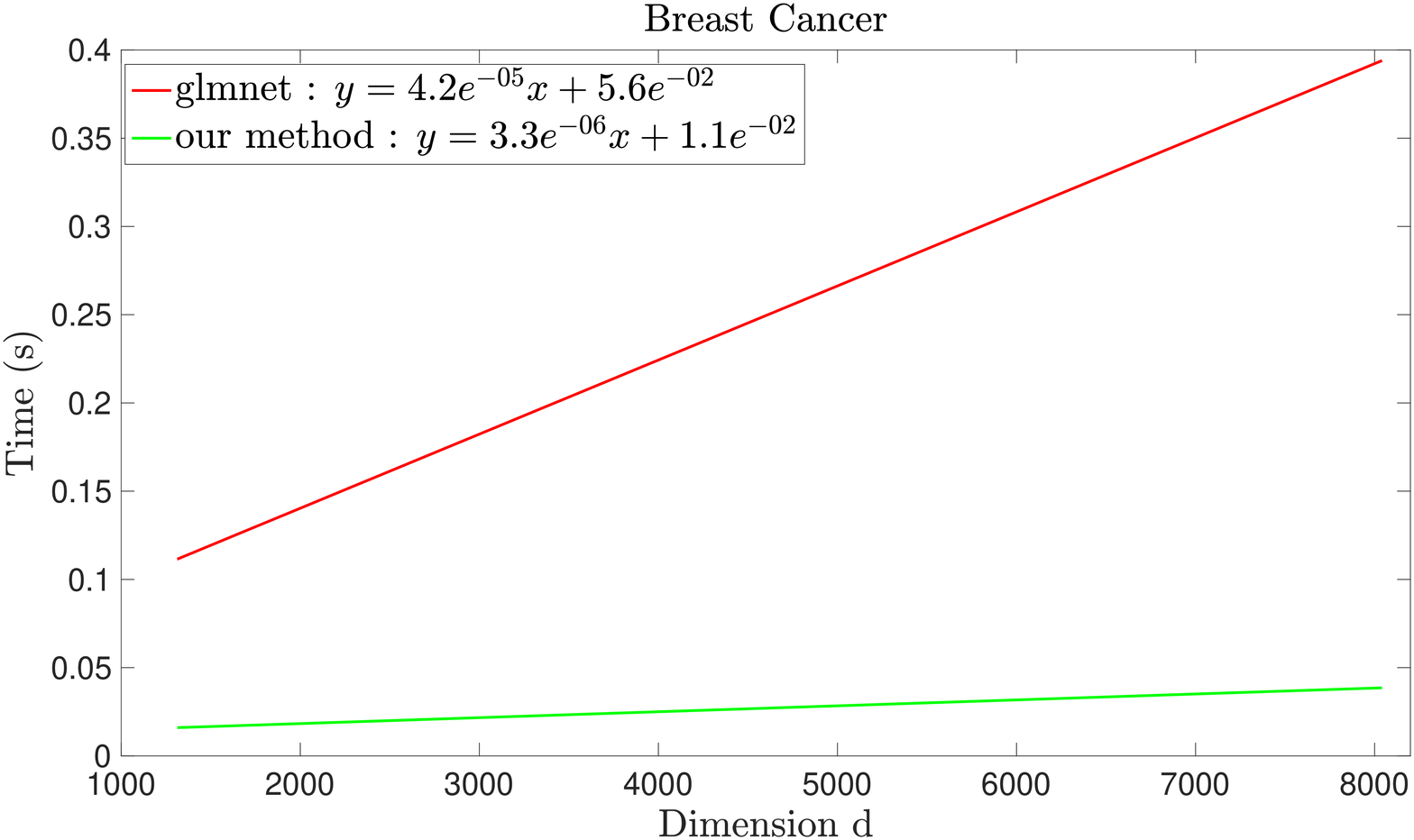}
\end{center}
\caption{Computing time as a function of the dimension.}
\label{fig:computingtime}
\end{figure}

Our classification implementation uses the logistic loss. 
Let $\|w\|_1\leq\eta$ be the surrogate sparsity
constraint. Fig.~\ref{fig:sparse} shows for different
values of the feature dimension that the number 
of nonzero coefficients 
$\|w\|_0$ decreases monotonically with the number of iterations.
Consequently, the sweep search over $\eta$ consists in stopping the 
iterations of the algorithm when $\|w\|_0 $ reaches value specified
a priori.
Our direct constrained strategy does not require the often
heuristic search for meaningful and interpretable Lagrange
multipliers. Moreover, we can improve 
processing time using sparse computing. Namely, at each iteration we
compute scalar products using only the sparse sub-vector of
nonzero values. We compare time processing using the breast
cancer data set ($n=295$ samples, $d=3022$) described in
Section~\ref{subsec:breast}. We provide time comparison using a
2.5 GHz Macbook Pro with an i7 processor and Matlab software. 
We report time processing in Table~\ref{tab:TimeBreast} using
glmnet software \cite{glmnet} and our method using either Matlab
$\ell^1$ or a standard mex $\ell^1$ file. Moreover, since the
vector $w$ is sparse, we provide mex-sparse and matlb-sparse
times using sparse computing.
\begin{table}
\begin{center}
\caption{Time comparison (Matlab and mex) versus glmnet \cite{glmnet}.}
\label{tab:TimeBreast}
\begin{tabular}{|c|c|c|c|c|c|c|}
\hline
&mex-sparse & mex & Matlab  &  Matlab-sparse & mex  \cite{glmnet}\\
\hline
Time(s)&0.0230&0.0559& 0.169 & 0.0729&0.198\\
\hline
\end{tabular}
\end{center}
\end{table}
Fig.~\ref{fig:computingtime} shows that our constrained method is
ten times faster than glmnet \cite{glmnet}. A numerical experiment 
is available in \cite{Barl17}.

A potentially competitive alternative $\ell^1$ constrained
optimization algorithms for solving the projection part of our
constrained classification splitting algorithm is that described 
in \cite{condat}. We plug the specific projection onto
the ball algorithm into our splitting algorithm. We  provide time
comparison (in seconds) in Table~\ref{time} for 
classification for the breast cancer data set ($d=3022$) 
described in Section~\ref{subsec:breast}.
\begin{table}
\begin{center}
\caption{Time comparison(s) with projection onto the ball
\cite{condat} for dimension $d=3022$ using Matlab.}
\label{time}
\begin{tabular}{|c|c|c|c|}
\hline
& Matlab &  Matlab sparse & ball \cite{condat} \\
\hline
Time (s) & 0.169&0.0729  &0.149 \\
\hline
\end{tabular}
\end{center}
\end{table}
Note that the most expensive part of our algorithm in terms of
computation is the evaluation of the gradient.
Although the projection onto the ball
\cite{condat} is faster than our projection, our method is
basically $12\%$ slower than the specific $\varphi_1$ constrained
method for dimension $d=3022$.  However, our sparse
implementation of scalar products is twice as fast. Moreover,
since the complexity of our method relies on the computation of
scalar products, it can be easily speed up using multicore CPU
or Graphics Processing Unit (GPU) devices, while the speed up of
the  projection on the ball  \cite{condat} using CPU or GPU is
currently  an open issue. 
In addition our method is more flexible since it can take
into account more sophisticated constraints such as $\varphi_2$,
$\varphi_3$, or any convex constraint.
We evaluate classification performance using area under the 
ROC curve (AUC). The result of Table~\ref{tab:Breast} show that 
our $\varphi_1$ constrained method outperforms
the $\varphi_1$ penalty method by $5.8\%$. 
Our $\varphi_2$ constraint improves slightly the AUC by
$1\%$ over the $\varphi_1$ constrained method. 
We also observe a significant improvement of our constrained 
$\varphi_2$ method over the penalty group Lasso
approach discussed in \cite{jpvert}.  
In addition, the main benefit of the $\varphi_2$ constraint is
to provide a set of connected genes which is more relevant
for biological analysis than the individual genes selected by 
the $\varphi_1$ constraint. 

\begin{table}
\begin{center}
\caption{Breast cancer AUC comparisons.}
\label{tab:Breast}
\begin{tabular}{|c|c|c|c|c|c|c|}
\hline
&glmnet \cite{glmnet}&Group Lasso
\cite{jpvert}&$\varphi_1$ & $\varphi_2$ \\
\hline
AUC ($\% $) &64.5& 66.7 & 71.3 & 72.3  \\
\hline
\end{tabular}
\end{center}
\end{table}

\begin{figure*}[t]
\begin{center}
\subfloat[]{
  \includegraphics[width=0.5\linewidth]{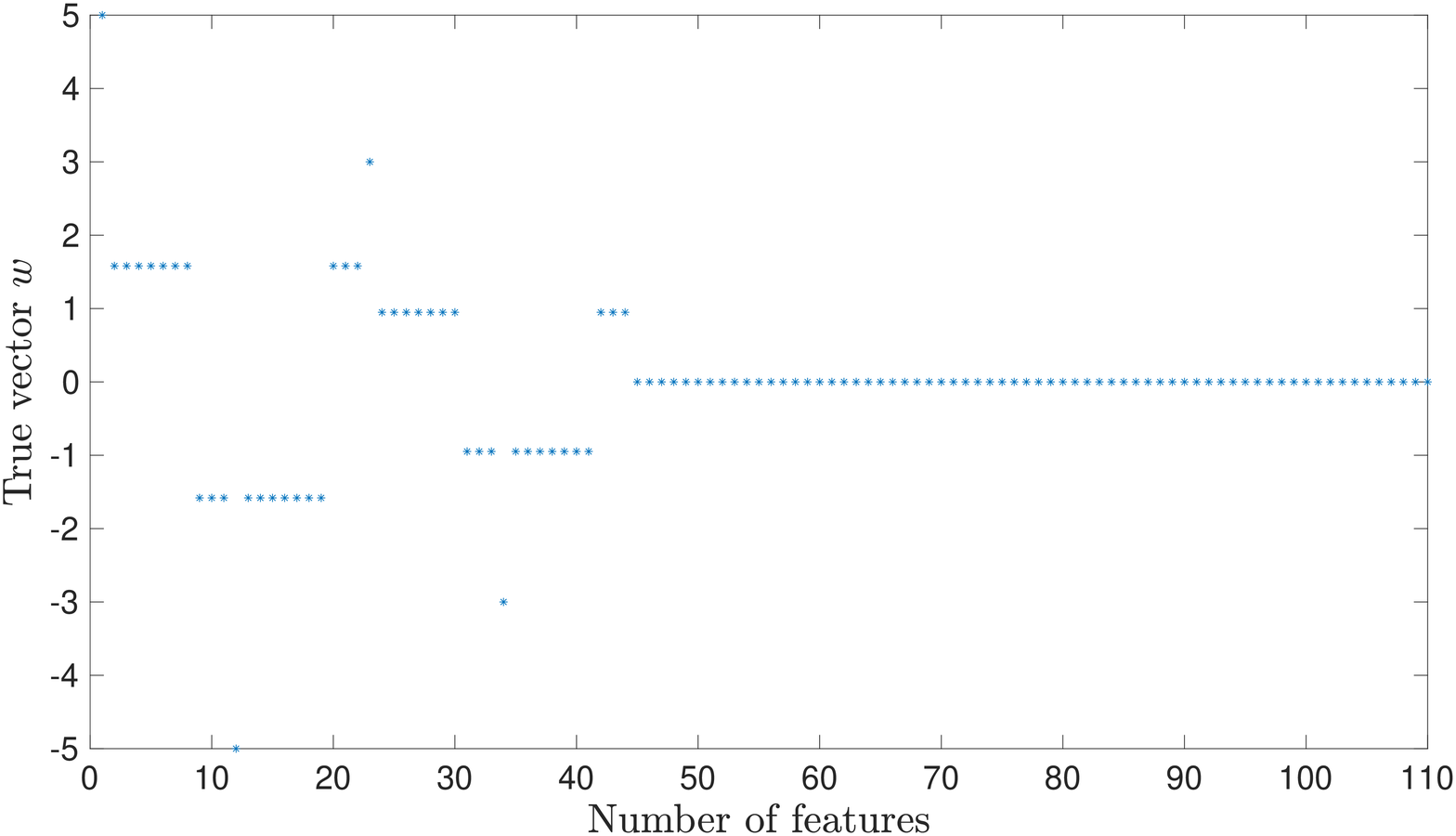}
  \label{subfig:estw1}
                     }
\subfloat[]{
  \includegraphics[width=0.5\linewidth]{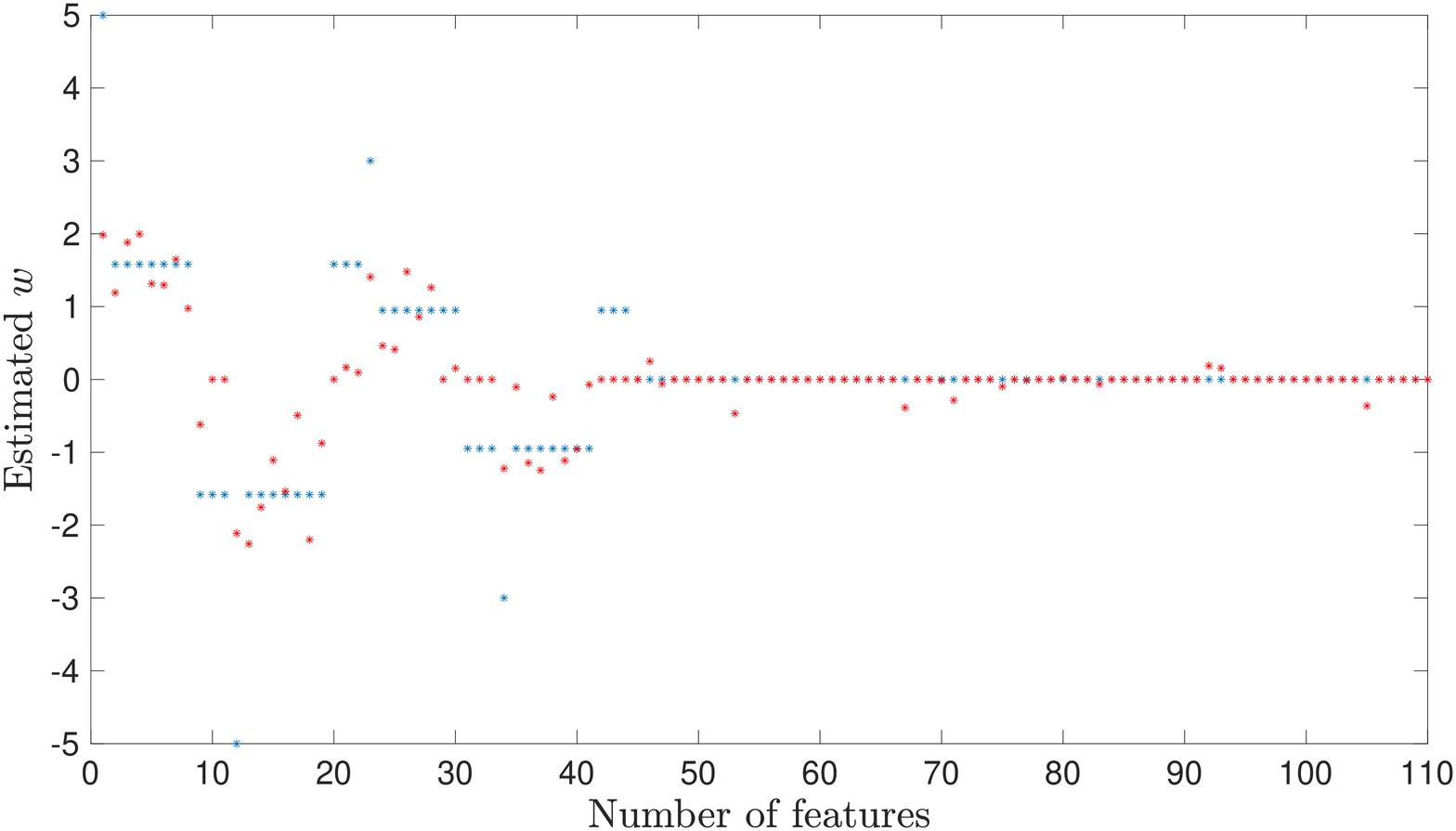}
  \label{subfig:estw2}
                     }
                     \\
\subfloat[]{
  \includegraphics[width=0.5\linewidth]{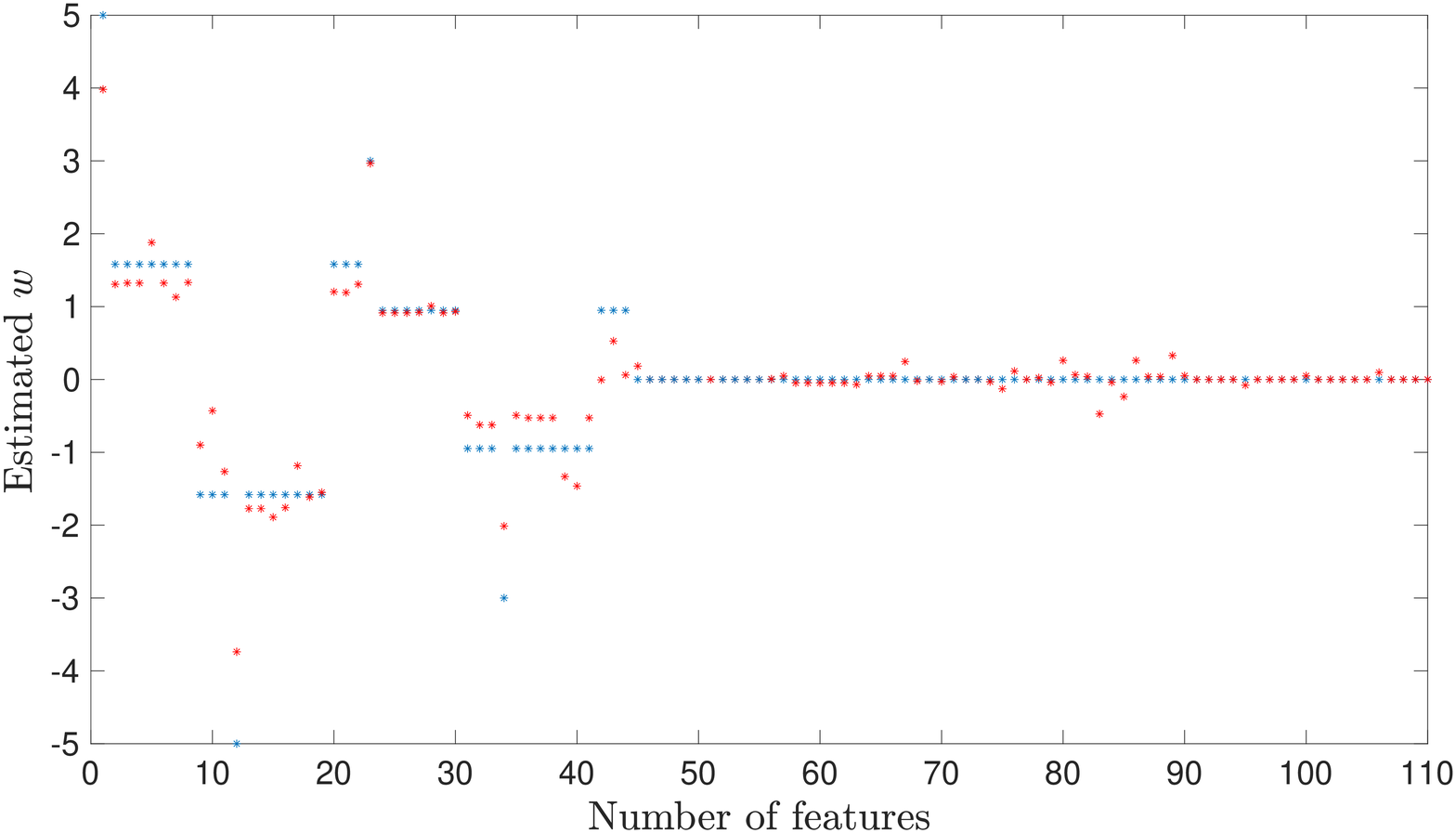}
  \label{subfig:estwR1}
                     }
\subfloat[]{
  \includegraphics[width=0.5\linewidth]{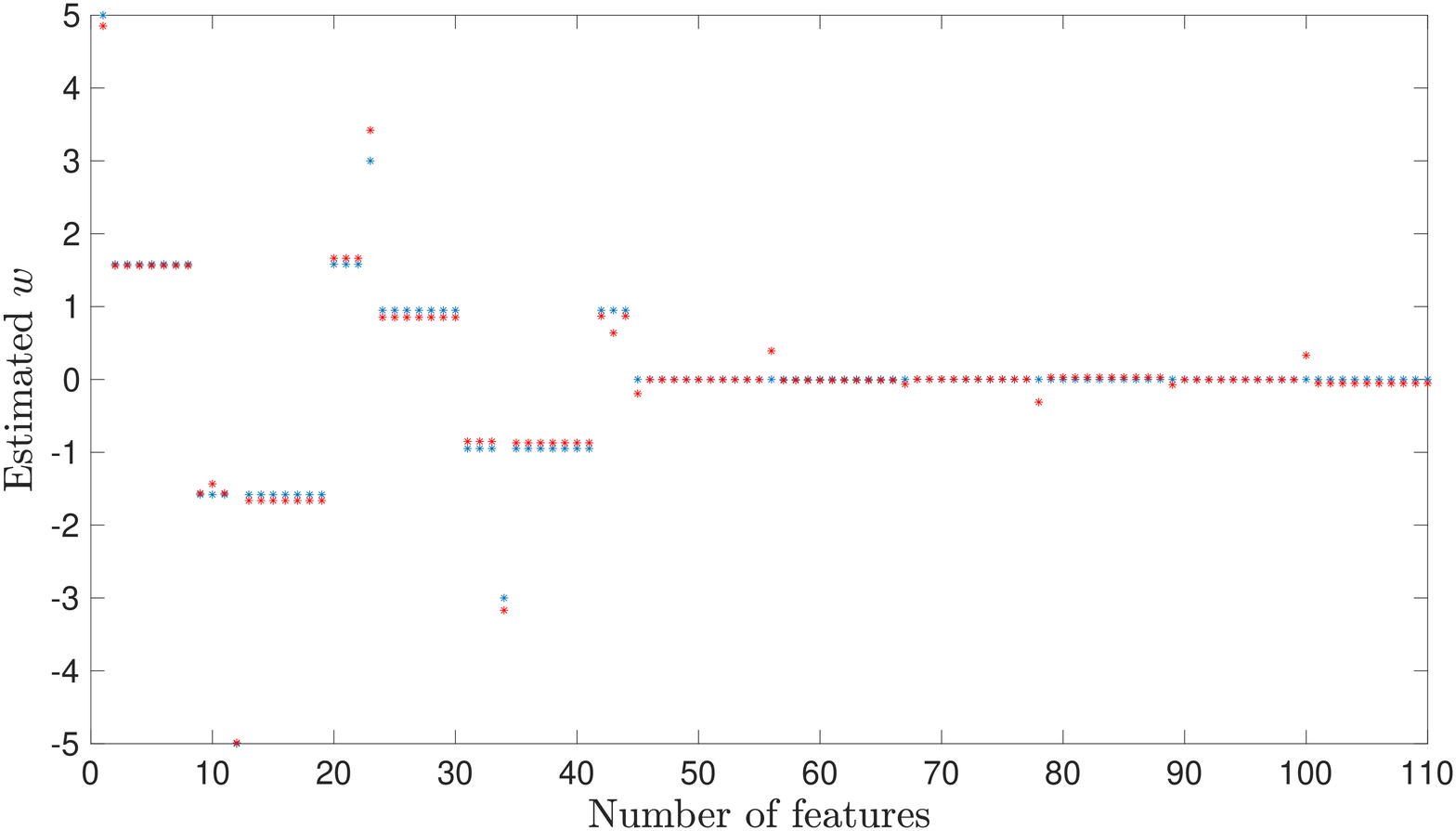}
  \label{subfig:estwR2}
                     }
\caption{Example~\ref{ex:3}: 
(a): True vector $w$.
(b): Estimation with the $\varphi_1$ constraint.
(c): Estimation with the $\varphi_2$ constraint.
(d): Estimation with the $\varphi_4$ constraint.}
\label{fig:estw}
\end{center}
\end{figure*}

\subsection{Comparison of various constraints for regression}
In biological applications, gene activation or inhibition are
well known and summarized in the 
ingenuity pathway analysis (IPA) database \cite{IPA}.
We introduce this biological a priori knowledge by replacing the
$\varphi_3$ constraint by
\begin{equation}
\label{e:def-varphi-bio}
\varphi_4\colon w\mapsto
\sum_{(i,j)\in\mathsf{S}}{|\omega_i-a_{i j} \omega_j|},
\end{equation}
where $a_{i j}=1$ if genes $i$ and $j$ are both activated or
inhibited, and $a_{i j}=-1$ if gene $i$ is activated and gene 
$j$ inhibited.
We compare the estimation of $w$ for Example~\ref{ex:3} using
$\varphi_1$ versus the $\varphi_2$ and $\varphi_4$ constraint. For
each fold, we estimate the regression vector $w$ on $100$ training
samples. Then we evaluate on new $100$ testing samples. We
evaluate regression using the mean square error (MSE) in the
training set and the predictive mean square error (PMSE) in the
test set. We use randomly half of the data for training and half
for testing, and then we average the accuracy over 50 random
folds.

\begin{figure}
\begin{center}
\includegraphics[width=\linewidth]{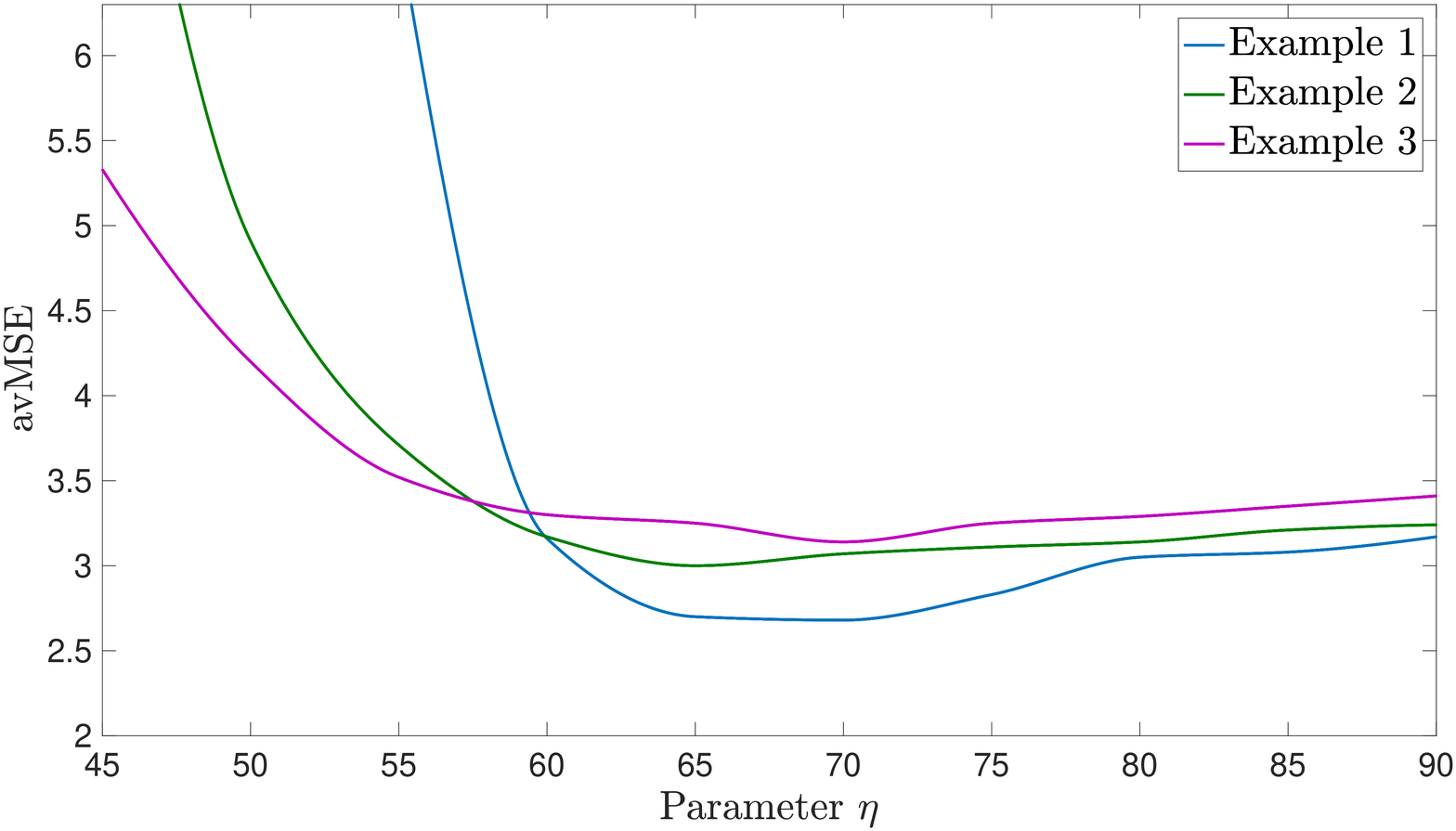}
\caption{$\varphi_1$ constraint for Examples~\ref{ex:1},
\ref{ex:2}, and \ref{ex:3}. Mean square error as a function of
the parameter $\eta$.}
\label{fig:Lasso}
\end{center}
\end{figure}

\begin{figure}
\begin{center}
\includegraphics[width=\linewidth]{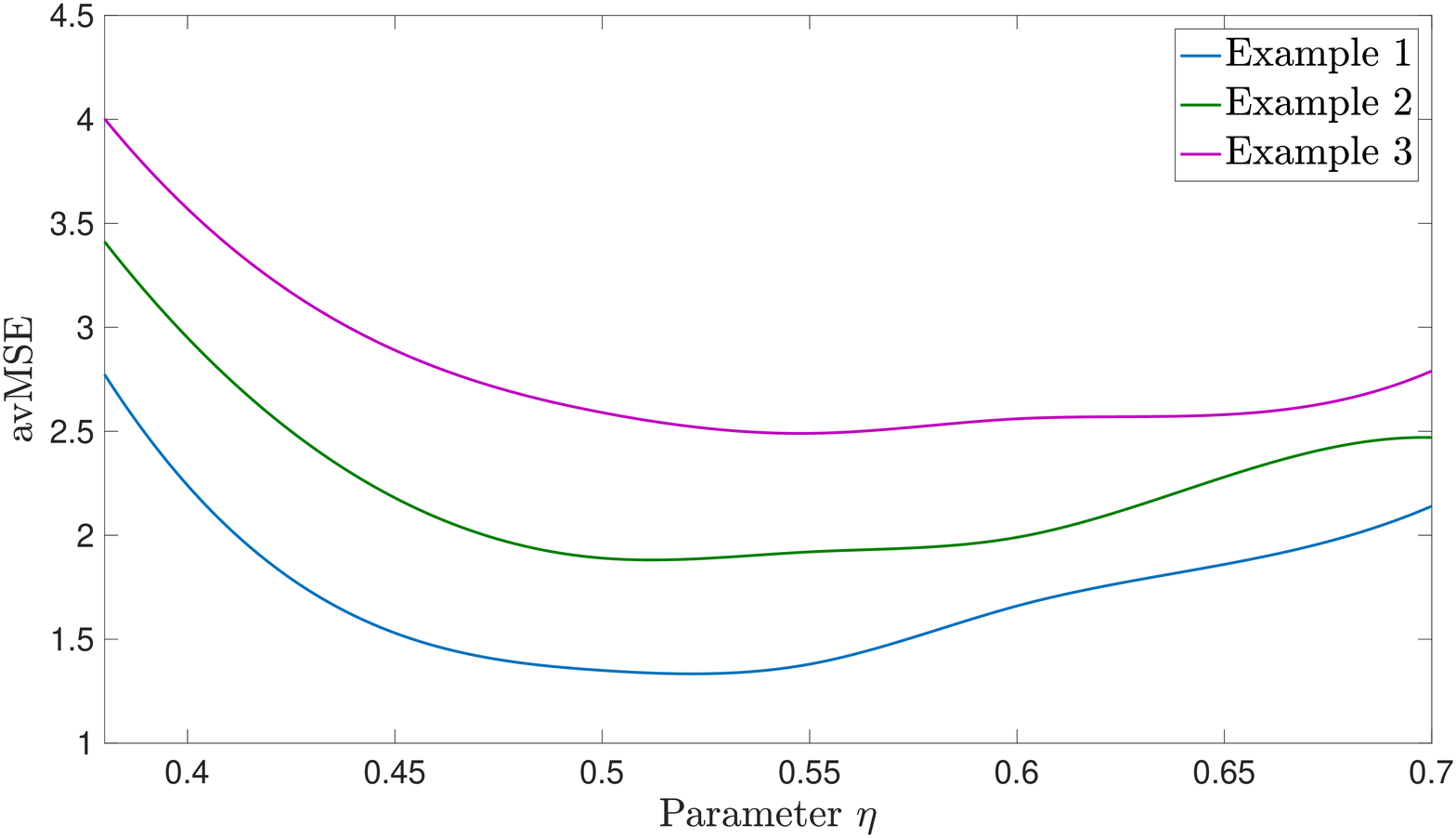}
\caption{$\varphi_2$ constraint for Examples~\ref{ex:1},
\ref{ex:2}, and \ref{ex:3}. Mean square error as a function of
the parameter $\eta$.}
\label{fig:graph}
\end{center}
\end{figure}

\begin{figure}
\begin{center}
\includegraphics[width=\linewidth]{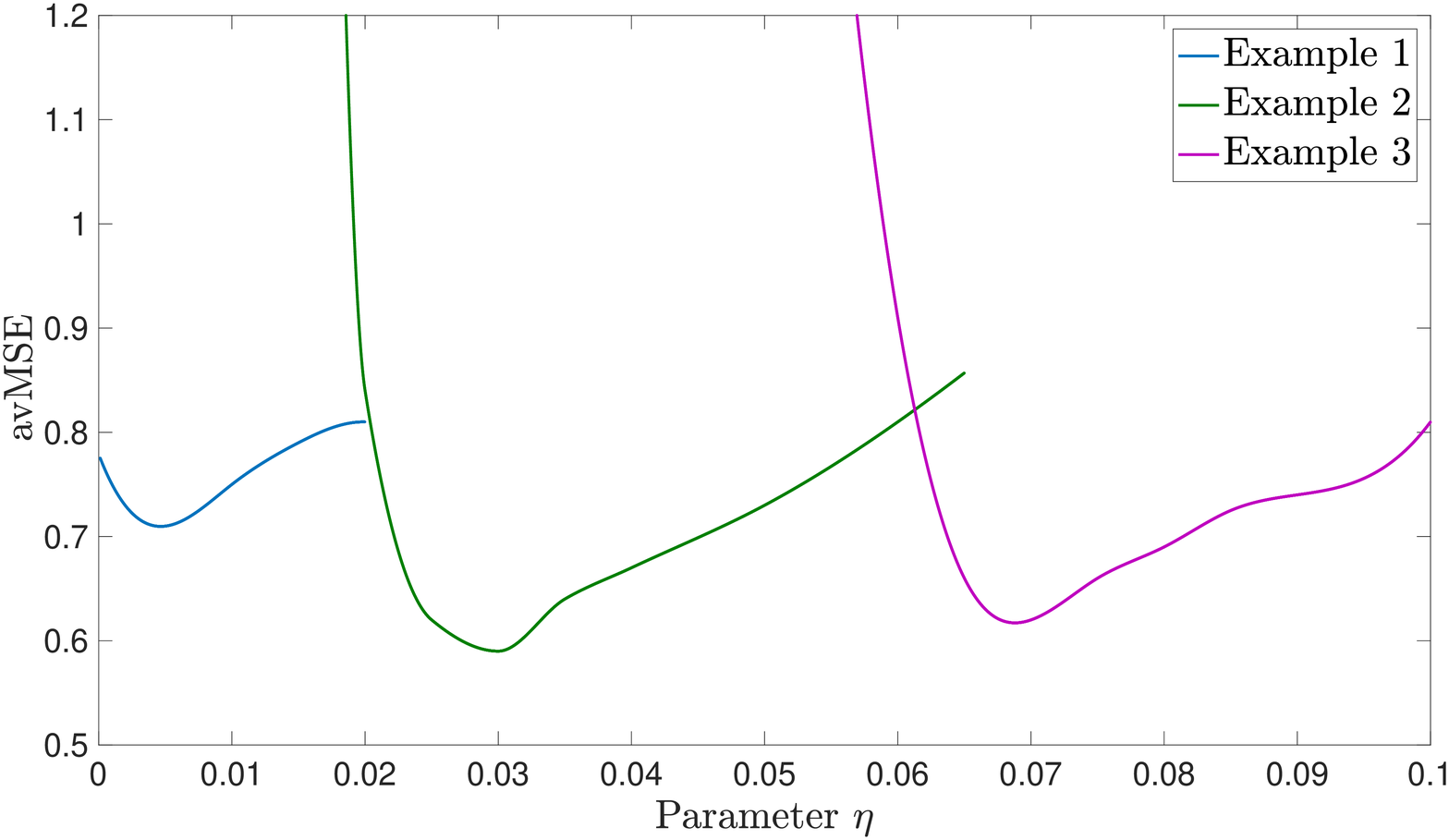}
\caption{$\varphi_4$ constraint for Examples~\ref{ex:1},
\ref{ex:2}, and \ref{ex:3}. Mean square error as a function of
the parameter $\eta$.}
\label{fig:dag}
\end{center}
\end{figure}

\begin{figure}
\begin{center}
\includegraphics[width=\linewidth]{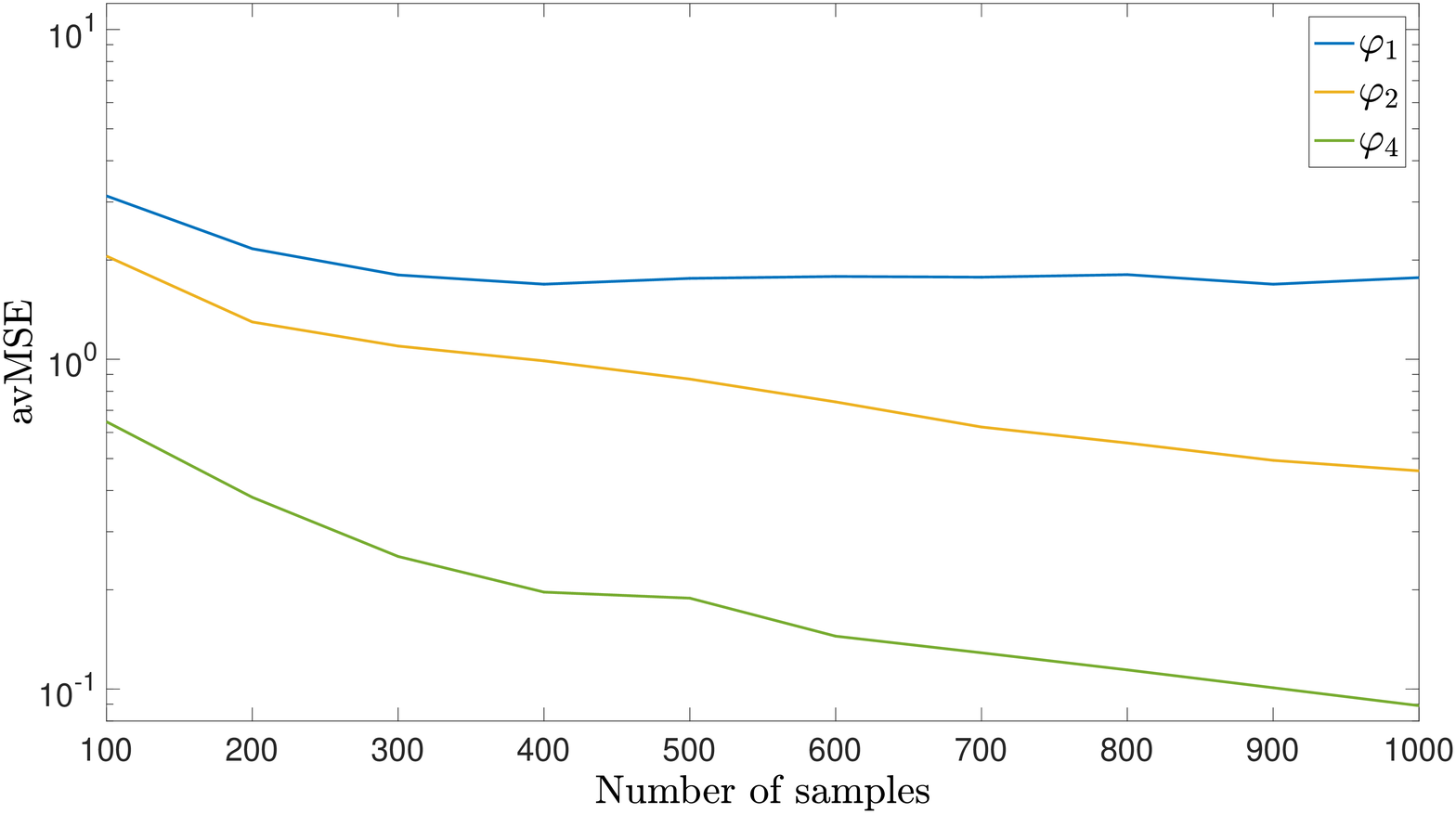}
\caption{MSE as a function of the number of samples $m$ for 
Example~\ref{ex:2}.}
\label{fig:MSE}
\end{center}
\end{figure}

We show in Fig.~\ref{subfig:estw1} the true regression vector
and, in Fig.~\ref{subfig:estw2}, the estimation 
using the $\varphi_1$ constraint for Example~\ref{ex:3}.  In
Fig.~\ref{subfig:estwR1} we show the results of the estimation
with the $\varphi_2$
constraint, and in Fig.~\ref{subfig:estwR2} with the $\varphi_4$
constraint. We provide for the three
examples the mean square error as a function of 
$\eta$ for $\varphi_1$ (Fig.~\ref{fig:Lasso}), $\varphi_2$
(Fig.~\ref{fig:graph}), and $\varphi_4$ (Fig.~\ref{fig:dag}). 
We report for Example~\ref{ex:2} in Fig.~\ref{fig:MSE} the 
estimation of the mean square error in the training 
set as a function of the number of training samples for the
$\varphi_1$, $\varphi_2$, and $\varphi_4$ constraint.
The $\varphi_4$ constraint outperforms both the $\varphi_2$
and the $\varphi_1$ constrained method. 
However, the selection of the parameter $\eta$ for 
constraint $\varphi_4$ is more challenging.

\section{Conclusion}
We have used constrained optimization approaches to promote
sparsity and feature grouping in classification and regression
problems. To solve these problems, we have proposed a new
efficient algorithm which alternates a gradient step on the data
fidelity term and an approximate projection step onto the
constraint set. We have
also discussed the generalization to multiple constraints.
Experiments on both synthetic and biological data show that 
our constrained approach outperforms penalty methods. Moreover, the
formulation using the $\varphi_4$ constraint outperforms
those using the pairwise $\varphi_2$ and the $\varphi_1$ 
constraint. 

\setcounter{equation}{0}
\renewcommand{\theequation}{A\arabic{equation}}
\section*{Appendix A -- The case of multiple constraints}
\label{A}
Let $\Phi$ be as in Problem~\ref{prob:2} and,
for every $j\in\{1,\ldots,p\}$, let $\varphi_j\colon\RR^d\to\RR$ 
be convex, let $\eta_j\in\RR$, and let
$\omega_j\in\zerounr$ be such that $\sum_{j=1}^p\omega_j=1$. 
Consider the problem 
\begin{equation}
\label{e:4}
\minimize{\substack{
\varphi_1(w)\leq\eta_1\\[-1mm]
~~~~\vdots\\[1mm]
\varphi_p(w)\leq\eta_p}}{\Phi(w)}.
\end{equation}
In other words,
$C=\bigcap_{j=1}^p\menge{w\in\HH}{\varphi_j(w)\leq\eta_j}$ in 
\eqref{e:prob2}. Let $k\in\NN$. 
For every $j\in\{1,\ldots,p\}$, let 
$s_{j,k}\in\partial\varphi_j(p_k)$ and set 
\begin{equation}
\label{e:e-radigue0}
p_{j,k}=
\begin{cases}
p_k+\displaystyle{\frac{\eta_j-\varphi_j(p_k)}{\|s_{j,k}\|^2}}s_{j,k}
&\text{if}\;\;\varphi_j(p_k)>\eta_j\\
p_k&\text{if}\;\;\varphi_j(p_k)\leq\eta_j.
\end{cases}
\end{equation}
Now define
\begin{equation}
\label{e:e-radigue1}
p_{k+1/2}=p_k+L_k
\Bigg(\Sum_{j=1}^p\omega_jp_{j,k}-p_k \Bigg),
\end{equation}
where
\begin{equation}
\label{e:e-radigue2}
L_k=\frac{\Sum_{j=1}^p\omega_j\|p_{j,k}-p_k\|^2}
{\Bigg\|\Sum_{j=1}^p\omega_jp_{j,k}-p_k\Bigg\|^2}.
\end{equation}
Then $p_{k+1}=Q(p_0,p_k,p_{k+1/2})\to P_C(p_0)$ 
\cite[Theorem~6.4]{Sico00} and therefore the generalization 
\begin{equation}
\label{jK$h18ghGa03-19x}
\begin{array}{l}
\text{for}\;n=0,1,\ldots\\
\left\lfloor
\begin{array}{l}
v_n=w_n-\gamma_n\nabla\Phi(w_n)\\
p_0=v_n\\
\begin{array}{l}
\hskip -2mm\text{for}\;k=0,1,\ldots,K_n-1\\
\left\lfloor
\begin{array}{l}
\zeta_k=\min_{1\leq j\leq p}\big(\eta_j-\varphi_j(p_k)\big)\\
\text{if}\;\zeta_k\geq 0\\
\begin{array}{l}
\left\lfloor
\text{terminate.}
\right.\\
\end{array}\\
\text{for}\;j=1,\ldots,p\\
\begin{array}{l}
\left\lfloor
s_{j,k}\in\partial\varphi_j(p_k)
\right.\\
\end{array}\\
\text{compute}\:p_{k+1/2}\;
\text{as in \eqref{e:e-radigue0}--\eqref{e:e-radigue2}}
\end{array}
\right.\\[2mm]
\end{array}\\
w_{n+1}=p_{K_n}.
\end{array}
\right.\\[2mm]
\end{array}
\end{equation}
of \eqref{wfwEkj7yUj01-23x} solves \eqref{e:4}.

\end{document}